\documentclass[microtype]{gtpart}

\usepackage{amsmath}
\usepackage{amsthm}
\usepackage{amsopn}
\usepackage{amssymb}
\usepackage[pdftex]{graphicx}
\usepackage[utf8]{inputenc}
\usepackage[T1]{fontenc}

\title{On the Slice Spectral Sequence}
\author{John Ullman}
\givenname{John}
\surname{Ullman}
\address{Department of Mathematics\\Massachusetts Institute of Technology\\\newline Cambridge, MA 02139\\USA}
\email{jrullman@math.mit.edu}
\urladdr{http://math.mit.edu/~jrullman/}


\arxivreference{arXiv:1206.0058 [math.AT]}
\arxivpassword{8w924}

\volumenumber{13}
\issuenumber{3}
\publicationyear{2013}
\papernumber{}
\startpage{1743}
\endpage{1755}
\doi{10.2140/agt.2013.13.1743}
\MR{}
\Zbl{}
\received{12 June 2012}
\revised{29 October 2012}
\accepted{12 November 2012}
\published{18 May 2013}
\publishedonline{}
\proposed{}
\seconded{}
\corresponding{}
\editor{}
\version{}

\newcommand{\ZZ}{\mathbb{Z}}
\newcommand{\RR}{\mathbb{R}}
\newcommand{\FF}{\mathcal{F}}
\newcommand{\PP}{\bar{P}}
\newcommand{\tr}{\bar{\tau}}

\newtheorem{thm}[equation]{Theorem}
\newtheorem{lem}[equation]{Lemma}
\newtheorem{prop}[equation]{Proposition}
\newtheorem{cor}[equation]{Corollary}
\newtheorem*{rem}{Remark}
\newtheorem*{warn}{Warning}

\numberwithin{equation}{section}

\newlength\tindent
\renewcommand{\indent}{\hspace*{\tindent}}

\begin{document}

\begin{abstract}
We introduce a variant of the slice spectral sequence which uses only regular slice cells, and state the precise relationship between the two spectral sequences. We analyze how the slice filtration of an equivariant spectrum that is concentrated over a normal subgroup is related to the slice filtration of its geometric fixed points, and use this to prove a conjecture of Hill on the slice filtration of an Eilenberg MacLane spectrum (\cite{Primer}). We also show how the (co)connectivity of a spectrum results in the (co)connectivity of its slice tower, demonstrating the "efficiency" of the slice spectral sequence.
\end{abstract}

\maketitle

\section{Introduction}\label{sec:intro}

\indent The slice spectral sequence in equivariant stable homotopy theory is analogous to the slice spectral sequence in motivic homotopy theory (see \cite{Voev}) and was an important tool in the recent solution of the Kervaire invariant problem (see \cite{HHR}). It is a spectral sequence of Mackey functors which converges strongly to the equivariant stable homotopy groups of the given spectrum. It is defined as the homotopy spectral sequence of the \emph{slice tower}, which is an analogue of the Postnikov tower. The slice tower is constructed using the apparatus of \emph{localizing subcategories}, which we review in Section~\ref{sec:loc}.\\
\indent In Section~\ref{sec:slice} we introduce the \emph{slice cells}, which play the same role in the construction of the slice tower that the sphere spectra play in the construction of the Postnikov tower, and define the \emph{regular slice tower} using the \emph{regular slice cells}. The precise relationship between the slice and regular slice towers is stated there as Proposition~\ref{prop:regslice}. We defer to~\cite{HHR} and~\cite{Primer} for basic facts and arguments about the slice filtration.\\
\indent In Sections~\ref{sec:Floc} and~\ref{sec:EMslice} we prove the conjectures of Hill (see~\cite{Primer}) concerning the slice towers of spectra which are \emph{concentrated over} a normal subgroup and the slice towers of Eilenberg MacLane spectra. In Section~\ref{sec:eff} we compute the slice towers of Eilenberg MacLane spectra in degree $-2$ and obtain bounds on the (co)connectivity of the slice tower of a spectrum with a given (co)connectivity. This demonstrates the relative "efficiency" of the slice spectral sequence.\\
\indent In what follows we implicitly use spectra in the sense of Lewis-May, where everything is fibrant, and the fibrations are the space-wise Serre fibrations (see \cite{LMS}). We work with finite groups $G$, $H$, etc.; restriction to a subgroup $H$ is denoted by $i^{*} _H$. The homotopy category of genuine $G$-spectra (indexed on a complete $G$-universe) is denoted by $Sp_G$, and the real regular representation of $G$ is denoted by $\rho_G$. We use $S^V$ to denote the one-point compactification of the representation $V$, and indicate Mackey functors with an underline (e.g. $\underline{M}$).

\section{Localizing Subcategories}\label{sec:loc}

\noindent We will say that a (nonempty) full subcategory $\tau$ of $Sp_{G}$ is \emph{localizing} if
\begin{itemize}
\item a spectrum isomorphic to an object of $\tau$ is in $\tau$,
\item $\tau$ is closed under taking cofibers and extensions,
\item $\tau$ is closed under wedge sums,
\item $\tau$ is closed under retract, and 
\item $\tau$ is closed under well-ordered homotopy colimits.
\end{itemize}
If $\tau$ is localizing then it contains the trivial spectrum and is closed under suspension and arbitrary homotopy colimits. From this it follows that the last two conditions above are redundant. See \cite{DF} for more. We will assume from now on that $\tau$ is generated by a set $T$ of spectra.\\
\indent Define $\tau\perp$ to be the full subcategory of spectra $X$ such that $[Y,X]=0$ for all $Y\in\tau$ (or equivalently, for $Y$ an iterated suspension of a generator). We write $\tau\perp X$ in place of $X\in\tau\perp$. Note that
\begin{itemize}
\item a spectrum isomorphic to an object of $\tau\perp$ is in $\tau\perp$,
\item $\tau\perp$ is closed under desuspension and taking fibers and extensions,
\item $\tau\perp$ is closed under products, and contains the trivial spectrum,
\item $\tau\perp$ is closed under retract and well-ordered homotopy limits, and
\item if $\tau$ has a set of compact generators then $\tau\perp$ is closed under wedge sums and directed homotopy colimits.
\end{itemize}

\noindent Recall from \cite{DF} that the inclusion $\tau\perp\subseteq Sp_{G}$ has a left adjoint $P^{\tau\perp}$. The construction is the familiar one of iteratively attaching nullhomotopies for maps from suspensions of generators. The construction must be iterated transfinitely in general, but may be iterated only countably many times if the generators are $\omega$-small. We have a functorial fiber sequence
\begin{align*}
	P_{\tau}X \to X \to P^{\tau\perp}X,
\end{align*}
letting $P_{\tau}X=Fib(X\to P^{\tau\perp}X)$, where the first map above is the terminal map to $X$ from a member of $\tau$ and the second is the initial map from $X$ to a member of $\tau\perp$. The spectrum $P_{\tau}X$ is also characterized by
\begin{itemize}
\item $P_{\tau}X \in \tau$, and 
\item $[\Sigma^{k} Y,P_{\tau}X]\to[\Sigma^{k} Y,X]$ is surjective for all $k\geq 0$ and injective for all $k\geq -1$ and all generators $Y$.
\end{itemize}
Thus we can construct $P_{\tau}X \xrightarrow{p} X$ as a colimit of stages $(P_{\tau} X)^{(i)} \xrightarrow{p_i} X$, as follows, when the generators are compact. First we let
\begin{align*}
	(P_{\tau} X)^{(0)} = \displaystyle \bigvee_{\substack{f : \Sigma^{k} Y \to X\\k\geq 0,Y\in T}} \textstyle \Sigma^k Y,
\end{align*}

\noindent where the sum runs over all maps $f : \Sigma^k Y \to X$ with $k \geq 0$ and $Y \in T$, and $p_0 = \vee_f f$. Next, assuming we have defined $(P_{\tau} X)^{(i)}$ and $p_i$, we define $(P_{\tau} X)^{(i+1)}$ as the cofiber of the map
\begin{align*}
	\displaystyle \bigvee_{\substack{f : \Sigma^k Y \to (P_{\tau} X)^{(i)}\\p_i f = 0,k \geq -1,Y\in T}} \textstyle \Sigma^k Y \xrightarrow{\vee_f f} (P_{\tau} X)^{(i)},
\end{align*}

\noindent where the sum runs over all maps $f : \Sigma^k Y \to (P_{\tau} X)^{(i)}$ with $k \geq -1$ and $Y \in T$ such that $p_i f = 0$, and $p_i$ is extended arbitrarily to $p_{i+1}$. In this construction, one is free to use only nonzero $f$ at each stage. In what follows, we will generally find the $P_{\tau}$ to be more useful than the $P^{\tau\perp}$.

\section {The Slice and Regular Slice Filtrations}\label{sec:slice}

\noindent We define the \emph{slice cells of dimension $k$} to be the spectra $G_{+}\wedge_{H}S^{n\rho_{H}}$ for $n|H|=k$ and $G_{+}\wedge_{H}S^{n\rho_{H}-1}$ for $n|H|-1=k$. Similarly, we define the \emph{regular slice cells} to be the slice cells of the first type listed above. Let $\tau_{k}$  (resp. $\tr_{k}$) be the localizing category generated by the slice cells (resp. regular slice cells) of dimension $\geq k$. We will write $\tau^G _n$, etc. if there is more than one group under consideration. The following facts are elementary:
\begin{itemize}
\item $\tau_{n}$ and $\tr_{n}$ are closed under induction and restriction, and thus under smashing with $(-1)$-connected spectra,
\item $\tau_{n} \subseteq \tau_{n-1}$, $\tr_{n} \subseteq \tr_{n-1}$,
\item $\tr_{n} \subseteq \tau_{n}, \tau_{n} \subseteq \tr_{n-(|G|-1)}$,
\item $S^{\rho_{G}} \wedge \tau_{n} \cong \tau_{n+|G|}, S^{\rho_{G}} \wedge \tr_{n} \cong \tr_{n+|G|}$.
\end{itemize}

\noindent See \cite{HHR} and \cite{Primer} for the basic arguments and results on the slice filtration. Less obvious is the following crucial fact.

\begin{prop}\label{prop:regslice}
For all $n$ we have $\Sigma\tau_{n} = \tr_{n+1}$.
\end{prop}
\begin{proof}
Since $\Sigma$ is an equivalence that maps cofiber sequences to cofiber sequences, it suffices to show that the suspension of a slice cell of dimension $k$ is in $\tr_{k+1}$ and that the desuspension of a regular slice cell of dimension $k$ is in $\tau_{k-1}$. The only nontrivial part is showing the inclusion $G_{+}\wedge_{H}S^{n\rho_{H}+1}\in\tr_{n|H|+1}$. We prove this by induction on $|G|$; the result is trivial for the trivial group. Thus we may assume the result for all proper subgroups of $G$. Since induction preserves the regular slice filtration, we may assume that $H=G$. Now take the cofiber sequence
\begin{align*}
	S(\rho_{G}-1)_{+} \to S^{0} \to S^{\rho_{G}-1}
\end{align*}
and smash with $S^{n\rho_{G}+1}$ to obtain
\begin{align*}
	S(\rho_{G}-1)_{+}\wedge S^{n\rho_{G}+1} \to S^{n\rho_{G}+1} \to S^{(n+1)\rho_{G}}.
\end{align*}
The spectrum on the left is built out of induced cells, and so the induction hypothesis implies that it is in $\tr_{n|G|+1}$. The spectrum on the right is in $\tr_{(n+1)|G|} \subseteq \tr_{n|G|+1}$. Thus the middle spectrum is in $\tr_{n|G|+1}$, as required.
\end{proof}

\begin{rem}
Amusingly, this implies that $\tau_n$ is generated by the \emph{irregular slice cells} of dimension $\geq n$.
\end{rem}

\begin{cor}\label{cor:sigmaconj}
For all $n$ we have natural isomorphisms
\begin{align*}
	\Sigma P_{\tau_{n}} \cong P_{\tr_{n+1}}\Sigma.
\end{align*}
\end{cor}

\begin{rem}
What the above implies is that the original slice construction is just a shifted version of the regular construction applied to the suspension. Thus we may view the regular construction as being fundamental, and the original construction as being derived from it.
\end{rem}

\begin{cor}\label{cor:susp}
For all $n$ we have inclusions
\begin{align*}
	\Sigma\tau_{n} \subseteq \tau_{n+1}, \Sigma\tr_{n} \subseteq \tr_{n+1}.
\end{align*}
\end{cor}

\noindent The proofs are immediate. Next we note that under certain circumstances, the slice and regular slice filtrations of a spectrum coincide.

\begin{prop}\label{prop:twosame}
If $P_{\tau_{n}}X \in \tr_{n}$ then $P_{\tau_{n}}X = P_{\tr_{n}}X$.
\end{prop}

\noindent The proof is immediate, considering the universal property which characterizes $P_{\tr_{n}}X$. We can also give a criterion in terms of the slices. We use the following notation for the $n$-slice:
\begin{align*}
	P^n_n X = Fib(P^{\tau_n\perp} X \to P^{\tau_{n-1}\perp} X) \cong Cofib(P_{\tau_{n+1}} X \to P_{\tau_n} X).
\end{align*}

\begin{prop}\label{prop:twosameslices}
If $P^k _k X \in \tr_{k}$ for all $k$ then $P_{\tau_{n}}X \in \tr_{n}$ for all $n$.
\end{prop}
\begin{proof}
Since $P_{\tau_{n+|G|-1}}X \in \tr_{n}$, the spectrum $P_{\tau_{n}}X$ has a finite filtration
\begin{align*}
	\ast \to P_{\tau_{n+|G|-1}}X \to ... \to P_{\tau_{n+1}}X \to P_{\tau_{n}}X
\end{align*}
such that the successive cofibers are in $\tr_{n}$.
\end{proof}

\noindent This brings us to the following important point.

\begin{rem}
The regular slice spectral sequence (or RSSS) is easily seen to have multiplicative pairings (one simply shows that $\tr_n \wedge \tr_m \subseteq \tr_{n+m}$). Also, the collection of regular slice cells is self-dual. Thus, we will observe below a kind of duality or symmetry about $0$ in the RSSS, while the slice spectral sequence (or SSS) is symmetric about $-1$. Furthermore, in favorable cases such as $K\RR$ (see \cite{KR}) and the various spectra constructed from $MU\RR$ in the solution of the Kervaire invariant problem (see \cite{HHR}), the two spectral sequences coincide. Thus, in many cases it may be more fruitful to work with the RSSS than the SSS. Note however that in cases where the two coincide, we are guaranteed different vanishing lines of slope $(|G|-1)$. For $t-s < 0$ the SSS guarantees a stronger vanishing line, while for $t-s > 0$ the RSSS does.
\end{rem}

\indent We can use Proposition~\ref{prop:regslice} to quickly derive a few more basic facts, using what is known about the $\tau_{n}$.

\begin{prop}\label{prop:facts}
The category $\tr_{k}$ consists of the $(k-1)$-connected spectra for $k=0,1$. The category $\tr_{2}$ consists of the connected spectra $X$ such that $[G_{+} \wedge S^{1}, X]=0$. For $k \geq 0$ and all $H \subseteq G$ we have $G/H_{+} \wedge S^{k} \in \tr_{k}$.
\end{prop}

\noindent The last statement is proved by induction on $k$, using Corollary~\ref{cor:susp} above. Finally, since we have $\tau_{0}=\tr_{0}$ and $S^{\rho_{G}} \wedge \tau_{n} \cong \tau_{n+|G|}, S^{\rho_{G}} \wedge \tr_{n} \cong \tr_{n+|G|}$, we conclude:

\begin{cor}\label{cor:periodicsame}
For all $n$ the categories $\tau_{n|G|}$ and $\tr_{n|G|}$ coincide.
\end{cor}

\indent From now on, we work with the RSSS instead of the SSS. Results for the SSS can be easily deduced from what follows by applying Proposition~\ref{prop:regslice}.\\
\indent If $X\in\tr_{n}$ we write $X \geq n$ or $X > n-1$, and if $\tr_{n}\perp X$ we write $X < n$ or $X \leq n-1$. We write $\PP_{n}$ in place of $P_{\tr_{n}}$ and $\PP^{n-1}$ in place of $P^{\tr_{n}\perp}$, so that we have functorial fiber sequences
\begin{align*}
	\PP_{n}X \to X \to \PP^{n-1}X.
\end{align*}

\noindent We denote the regular $n$-slice of $X$ by $\PP_n^n X$.

\section{Local Spectra}\label{sec:Floc}

\noindent Let $\FF$ be a family of subgroups of $G$, and $\FF'$ its complement in the family of all subgroups. As usual we denote by $E\FF$  the universal $\FF$-space and $\tilde{E}\FF$ its unreduced suspension. Recall that a spectrum $X$ is called \emph{$\FF'$-local} if one of the following equivalent conditions hold:
\begin{itemize}
\item $E\FF_{+} \wedge X \cong 0$,
\item $X \cong \tilde{E}\FF \wedge X$,
\item $i^{*}_{H} X \cong 0$ for all $H \in \FF$,
\item $\pi^{H}_{n}X=0$ for all $n$ and all $H \in \FF$.
\end{itemize}

\noindent The inclusion of the full subcategory of $\FF'$-local spectra, which we denote by $Sp^{\FF'} _G$, has a left adjoint, given by $\tilde{E}\FF \wedge (\cdot)$. We call this \emph{$\FF'$-localization.} Now let $\tau$ denote a localizing subcategory which is generated by a set $T$. We have the following general fact.

\begin{thm}\label{thm:Flocal}
If $\tau$ is closed under $\FF'$-localization, then for any $\FF'$-local spectrum $X$, $P_{\tau}X$ and $P^{\tau\perp}X$ are $\FF'$-local.
\end{thm}
\begin{proof}
The second statement follows from the first. Let $\tau^{\FF'}$ denote the localizing subcategory generated by $\tilde{E}\FF \wedge T$; by hypothesis, this is contained in $\tau$. Also, any element of $\tau^{\FF'}$ has a (possibly transfinite) filtration whose succesive cofibers are wedges of suspensions of elements of $\tilde{E}\FF \wedge T$, so $\tau^{\FF'}$ consists of $\FF'$-local spectra. This implies that $\tau^{\FF'} \subseteq \tilde{E}\FF \wedge \tau$. Conversely, suppose that $X$ has a (possibly transfinite) filtration whose succesive cofibers are wedges of suspensions of elements of $T$. Smashing this filtration with $\tilde{E}\FF$, we obtain a filtration for $\tilde{E}\FF \wedge X$ with succesive cofibers that are wedges of suspensions of elements of $\tilde{E}\FF \wedge T$. Thus we have $\tilde{E}\FF \wedge \tau \subseteq \tau^{\FF'}$ as well, so that $\tau^{\FF'} = \tilde{E}\FF \wedge \tau = Sp^{\FF'} _G \cap \tau$.\\
\indent Now let $X$ be $\FF'$-local, and consider the map $P_{\tau^{\FF'}}X \to X$. The spectrum $P_{\tau^{\FF'}}X$ is in $\tau$ by the above. Furthermore, for any $Y \in \tau$ the map
\begin{align*}
	[Y, P_{\tau^{\FF'}}X] \to [Y, X]
\end{align*}
is isomorphic to the map
\begin{align*}
	[\tilde{E}\FF \wedge Y, P_{\tau^{\FF'}}X] \to [\tilde{E}\FF \wedge Y, X]
\end{align*}
since both $X$ and $P_{\tau^{\FF'}}X$ are $\FF'$-local. The above map is an isomorphism; thus, $P_{\tau^{\FF'}}X \to X$ satisfies the required universal property. That is, $P_{\tau}X \cong P_{\tau^{\FF'}}X$, so it is $\FF'$-local.
\end{proof}

\noindent Since the categories $\tr_{n}$ satisfy the above criterion, we immediately get the following:

\begin{cor}\label{cor:sliceFlocal}
The RSSS for an $\FF'$-local spectrum is $\FF'$-local.
\end{cor}

\begin{rem}
This corollary can be proven much more simply by using the fact that the RSSS construction commutes with restriction functors; in fact, if $i^* _H X = 0$ then $i^* _H \PP_n X = 0$ for all $n$. However, we will need the more precise arguments given above in what follows.
\end{rem}
\begin{warn}
The example of $G=\ZZ/2\ZZ$, $\FF=\{e\}$, $X=K\RR$ shows that taking (regular) slices does not commute with localization in general. In fact, $\tilde{E}\FF \wedge X = 0$, but the localizations of the (regular) slices of $X$ are not zero (see \cite{KR}).
\end{warn}
\indent Now suppose that $N$ is a normal subgroup of $G$, and let $\FF[N]$ denote the family of subgroups that do not contain $N$. Recall that $Sp^{\FF[N]'} _G$ is the category of spectra whose homotopy groups are \emph{concentrated over $N$}, and that this is equivalent to the category of $G/N$-spectra (see \cite{LMS}). The equivalence is given by the $N$-fixed point functor, which is equal to the $N$-geometric fixed point functor $\Phi^{N}$ on $Sp^{\FF[N]'} _G$. Following Hill (\cite{Primer}), we call the inverse equivalence \emph{pullback}, and denote it by $\phi^{*}_{N}$. We recall the following basic fact, which may be proved by noting that $G \wedge_H (\cdot)$ and $\Phi^N$ are left adjoint functors whose right adjoints, $i^{*} _H$ and $\phi^{*} _N$, fit into a commutative square with $i^{*} _{H/N}$ and $\phi^{*} _N : Sp_{H/N} \to Sp_H$.

\begin{lem}\label{lem:geominduce}
There are natural isomorphisms as below.
\begin{align*}
	\Phi^{N}(G_{+} \wedge_{H} X) \cong \begin{cases}
	G/N_{+} \wedge_{H/N} \Phi^{N}(X) &\mbox{if } H \supseteq N \\
	0 &\mbox{if } H \not\supseteq N \end{cases}
\end{align*}
\end{lem}

\noindent In particular, the homotopy groups of $\Phi^{N}X$ are just the (possibly) nonzero homotopy groups of $X$, when $X$ is concentrated over $N$.  We can now prove the following results.

\begin{thm}\label{thm:GoverNslicefilt}
After saturating in isomorphism classes, we have
\begin{align*}
	\phi ^* _N (Sp_{G/N}) \cap \tr^G _m = \phi^* _N (\tr^{G/N} _{\lceil m/|N| \rceil}).
\end{align*}
\end{thm}
\begin{proof}
First, note that the above intersection is $Sp^{\FF[N]'} _G \cap \tr_m = \tr^{\FF[N]'}_{m}$. Now let $T_{m}$ denote the set of regular slice cells of dimension $\geq m$. Then $\tr^{\FF[N]'}_{m}$ is the localizing category generated by $\tilde{E}\FF[N] \wedge T_m$. Since $Sp^{\FF[N]'} _G$ is a triangulated subcategory of $Sp_{G}$ which is closed under wedge sums, we can regard $\tr^{\FF[N]'}_{m}$ as the localizing subcategory of $Sp^{\FF[N]'} _G$ generated by $\tilde{E}\FF[N] \wedge T_m$. Now, the geometric fixed point functor is an equivalence of triangulated categories, so we can immediately identify $\Phi^N (\tr^{\FF[N]'}_{m})$ as the localizing subcategory of $Sp_{G/N}$ generated by $\Phi^{N}(T_{m})$. To determine this category, we may begin by throwing out the elements of $\Phi^{N}(T_{m})$ that are trivial. By Lemma~\ref{lem:geominduce}, the nonzero elements are
\begin{align*}
	\Phi^{N}(G_{+} \wedge_{H} S^{k\rho_{H}}) \cong G/N_{+} \wedge_{H/N} S^{k\rho_{H/N}}
\end{align*}
for $H \supseteq N, k|H| \geq m$. These are all regular slice cells, and a regular slice cell of dimension $d=k|H/N|$ is in this set when $k|H|=d|N| \geq m$; that is, when $d \geq \lceil m/|N| \rceil$. Thus we have $\Phi^N (\tr^{\FF[N]'}_{m})=\tr_{\lceil m/|N| \rceil}$.
\end{proof}

\begin{cor}\label{cor:GoverNslices}
Let $X=\phi^* _N (Y)$ be concentrated over $N$. Then $\PP_m X = \phi^* _N (\PP_{\lceil m/|N| \rceil} Y)$. Thus the regular $k|N|$-slice of $X$ is the pullback of the regular $k$-slice of $Y$, and all other regular slices of $X$ are zero.
\end{cor}
\begin{proof}
We know that $\PP_m X$ is concentrated over $N$. Thus $\PP_m X \to X$ is terminal among maps to $X$ from spectra in $\phi ^* _N (Sp_{G/N}) \cap \tr^G _m = \phi^* _N (\tr^{G/N} _{\lceil m/|N| \rceil})$. The second statement follows from the first.
\end{proof}

\begin{rem}
The above result was originally conjectured by Mike Hill in~\cite{Primer}; it is stated there (in terms of the \emph{irregular} slice filtration) as Conjecture~4.11. Theorem~4.9 of that paper gives the correct upper bound for spectra that are pulled back from a quotient group, while Theorem~4.12 gives a non-optimal lower bound; the above corollary remedies this situation. Hill also gave proofs of the special cases where $N = G$ and where $[G : N] = 2$; see Theorem~6.14 and Corollary~4.14 of~\cite{Primer}, respectively.
\end{rem}

\noindent The following corollary will be useful for doing inductive proofs in the sequel.

\begin{cor}\label{cor:GoverGslicefilt}
If $X$ is concentrated over $G$ and $(n-1)$-connected then $X \geq n|G|$.
\end{cor}
\begin{proof}
$X^{G}$ is $(n-1)$-connected; hence, it is $\geq n$, so its lowest possible nonzero regular slice is the regular $n$-slice. By Theorem~\ref{thm:GoverNslicefilt}, the lowest possible nonzero regular slice of $X$ is the regular $n|G|$-slice.
\end{proof}

\section{The Slices of an Eilenberg MacLane Spectrum}\label{sec:EMslice}

\noindent Let $\underline{M}$ be a Mackey functor. Following Hill (see \cite{Primer}), we define a filtration of $\underline{M}$ by
\begin{align*}
	F^{k}\underline{M}(G/H) = \{x \in \underline{M}(G/H) : i^*_J x = 0\text{ in }\underline{M}(G/J)\text{ for }J \subseteq H, |J| < k\}.
\end{align*}
(Actually, this differs from Hill's filtration by a shift.) Evidently, $F^{k}\underline{M}(G/H)=0$ when $|H| < k$. To determine the regular slice tower for $\Sigma H\underline{M}$ we need two preliminary lemmas.

\begin{lem}\label{lem:EMgeq}
If $k \geq 1$ and $\underline{M}(G/H)=0$ for $|H| < k$ then $\Sigma H\underline{M} \geq k$.
\end{lem}
\begin{proof}
We proceed by induction on $|G|$; the result is trivial for the trivial group. We may also assume that $1 < k \leq |G|$, since $\underline{M}=0$ for $k > |G|$. Consider the cofiber sequence
\begin{align*}
	\PP_{k}(\Sigma H\underline{M}) \to \Sigma H\underline{M} \to C.
\end{align*}
The spectrum $C$ is evidently $(1-1)$-connected, and is concentrated over $G$ by the induction hypothesis. Therefore, by Corollary~\ref{cor:GoverGslicefilt}, we have ${C \geq |G| \geq k}$.
\end{proof}

\begin{lem}\label{lem:regcoh}
There is a natural isomorphism
\begin{align*}
	[G_+ \wedge_H S^{\rho_H}, \Sigma H\underline{M}] = \{x \in \underline{M}(G/H) : \text{all proper restrictions of }x\text{ are }0\}.
\end{align*}
\end{lem}
\begin{proof}
Firstly, the above group is isomorphic to $[S^{\rho_H - 1}, i^* _H H\underline{M}]$, which is isomorphic to $[(S^{\rho_H - 1})^{[1]}, H(i^* _H \underline{M})]$, where $(-)^{[1]}$ denotes taking the $1$-skeleton. There is a cofiber sequence
\begin{align*}
	(S(\rho_H-1))^{[0]} _+ \to S^0 \to (S^{\rho_H - 1})^{[1]} \to \Sigma (S(\rho_H-1))^{[0]} _+ .
\end{align*}
Now the space $S(\rho_H-1)$ can be given a cell structure with at least one zero-cell of the form $H/J$ for each proper subgroup $J$ of $H$, so we can take $(S(\rho_H-1))^{[0]} _+$ to be a wedge of copies of $H/J_+$, with at least one summand for each $J \subsetneq H$. The result follows by mapping the above cofiber sequence into $H(i^* _H \underline{M})$.
\end{proof}

\noindent We now prove the conjecture of Hill (\cite{Primer}) concerning the regular slice filtration for $\Sigma H\underline{M}$.

\begin{thm}\label{thm:EMplusslicetower}
The regular slice tower for $\Sigma H\underline{M}$ is given by 
\begin{align*}
	\PP_{k}(\Sigma H\underline{M}) = \Sigma H(F^k \underline{M}).
\end{align*}
\end{thm}
\begin{proof}
We have a cofiber sequence
\begin{align*}
	\Sigma H(F^k \underline{M}) \to \Sigma H\underline{M} \to \Sigma H(\underline{M}/F^k \underline{M}).
\end{align*}
By Lemma~\ref{lem:EMgeq}, the spectrum on the left is $\geq k$, so it suffices to show that the spectrum on the right is $< k$. Since it is an Eilenberg MacLane spectrum in dimension $1$, the only nontrivial part is to show that
\begin{align*}
	[G_+ \wedge_H S^{\rho_H}, \Sigma H(\underline{M}/F^k \underline{M})] = 0
\end{align*}
for all subgroups $H$ such that $|H| \geq k$. This group is identified by Lemma~\ref{lem:regcoh}. Hence, let $x \in \underline{M}/F^k \underline{M}(G/H)$ for some $H$ of order $\geq k$ and suppose that $i^* _J x = 0$ for all $J \subsetneq H$. Let $x$ be represented by $y \in \underline{M}(G/H)$. For any subgroup $J$ of $H$ of order $< k$, $J$ is proper, so $i^* _J x = 0$ in $\underline{M}/F^k \underline{M}(G/J) \cong \underline{M}(G/J)$. Thus, $i^* _J y = 0$ as well; that is, $y \in F^k \underline{M}(G/H)$, so $x=0$.
\end{proof}

\begin{rem}
The above result for cyclic groups of prime power order was previously given by Hill (see Theorem~5.8 of~\cite{Primer}). Hill was led to conjecture the above result by the observation that the filtration he defined on Mackey functors makes sense for arbitrary finite groups.
\end{rem}

\noindent We will compute the regular slice tower for a $\Sigma^{-1} H\underline{M}$ in the next section.

\section{Efficiency of the Slice Spectral Sequence}\label{sec:eff}

\noindent Let $n > 0$, and consider $G_+ \wedge_H S^{n\rho_H}$. This has a cell decomposition into cells of type $G/J$ and dimension $k$, where $k|J| \leq n|H|$. Taking the Spanier-Whitehead dual decomposition for $G_+ \wedge_H S^{-n\rho_H}$, we find that it is in the localizing subcategory generated by the spectra $G/J_+ \wedge S^{-k}$ with $k|J| \leq n|H|$. Furthermore, we have $G/J_+ \wedge S^{-k} \cong G_+ \wedge_J (S^{-k\rho_J} \wedge S^{k(\rho_J - 1)}) \geq -k|J|$. We obtain the following.

\begin{thm}\label{thm:spheregen}
Let $n \geq 0$. The category $\tr_{-n}$ is generated by the spectra $G/H_+ \wedge S^k$ for $k \geq 0$ and the spectra $G/H_+ \wedge S^{-k}$ for $k > 0, k|H| \leq n$. Also, $\tr_1$ is generated by the spectra $G/H_+ \wedge S^k$ for $k \geq 1$ and $H \subseteq G$.
\end{thm}

\begin{thm}\label{thm:conn}
Let $n \leq 1$, and suppose that $\underline{\pi}_i X = 0$ for $i < n$. Then this holds for all $\PP_k X$, $\PP^k X$, and $\PP^k _k X$.
\end{thm}
\begin{proof}
It suffices to prove that $\PP_m X$ is $(n-1)$-connected for all $m$. This holds trivially for $m > 0$, so we may assume that $m \leq 0$. We begin to construct $\PP_m X$ with a wedge of suspensions of the generators given in Theorem~\ref{thm:spheregen}, one for each (nonzero) map into X. Thus this first stage uses only spheres of dimension $\geq n$, so it is $(n-1)$-connected. Next, we iteratively kill the kernel on suspensions (by $j \geq -1$) of these generators, with one summand for each (nonzero) map in the kernel. Inductively, we easily see that at each stage we are only using spheres of dimension $\geq n$, so that each stage is $(n-1)$-connected.
\end{proof}
\noindent Note that, in the above construction, the map from $\underline{\pi}_n$ of each stage to $\underline{\pi}_n$ of the next stage is surjective. Hence, recalling that the regular slice filtration may be given by
\begin{align*}
	F^m \underline{\pi}_n X = im(\underline{\pi}_n \PP_m X \to \underline{\pi}_n X),
\end{align*}
we have the following corollary.

\begin{cor}\label{cor:htpyslicefilt}
Let $n > 0$ and suppose that $X$ is $(-n-1)$-connected. Then $F^{-m} \underline{\pi}_{-n} X$ is the sub-Mackey functor of $\underline{\pi}_{-n} X$ generated by the elements of $\underline{\pi}_{-n} X(G/H)$ for $|H| \leq m/n$.
\end{cor}

\noindent Next, for any Mackey functor $\underline{M}$ and any real number $c$, we define $F_c \underline{M}$ to be the sub-Mackey functor of $\underline{M}$ generated by the elements of $\underline{M}(G/H)$ for $|H| \leq c$.

\begin{thm}\label{thm:EMnegslicetower}
For any Mackey functor $\underline{M}$ and integer $n$ we have
\begin{gather*}
	\PP_n (\Sigma^{-1} H\underline{M}) = \Sigma^{-1} H(F_{-n} \underline{M}),\\
	\PP^n _n (\Sigma^{-1} H\underline{M}) = \Sigma^{-1} H(F_{-n} \underline{M} / F_{-n-1} \underline{M}).
\end{gather*}
\end{thm}
\begin{proof}
The second statement follows from the first. Now, $\Sigma^{-1} H\underline{M} < 0$, so both spectra are clearly zero for $n \geq 0$. Hence, let $n < 0$. Since the spectra $G/H_+ \wedge S^j$ are $> n$ for all $j \geq 0$, we have
\begin{gather*}
	\underline{\pi}_j \PP_n (\Sigma^{-1} H\underline{M}) \cong \underline{\pi}_j (\Sigma^{-1} H\underline{M}) = 0,\\
	\underline{\pi}_j \PP^{n-1} (\Sigma^{-1} H\underline{M}) = 0
\end{gather*}
for all $j \geq 0$. Furthermore, by Theorem~\ref{thm:conn} the spectra $\PP_n (\Sigma^{-1} H\underline{M})$ and $\PP^{n-1} (\Sigma^{-1} H\underline{M})$ are $(-2)$-connected; hence, they are Eilenberg MacLane spectra in dimension $-1$. This implies that the map
\begin{align*}
	\underline{\pi}_{-1} \PP_n (\Sigma^{-1} H\underline{M}) \to \underline{\pi}_{-1} (\Sigma^{-1} H\underline{M}) = \underline{M}
\end{align*}
is injective. Thus, by Corollary~\ref{cor:htpyslicefilt} applied with $n=1$ we have
\begin{align*}
	\underline{\pi}_{-1} \PP_m (\Sigma^{-1} H\underline{M}) \cong F^m \underline{M} = F_{-m} \underline{M}.
\end{align*}
\end{proof}

\noindent Analogous to Theorem~\ref{thm:conn} in this section, we have the following.

\begin{thm}\label{thm:coconn}
Let $n \geq -1$, and suppose that $\underline{\pi}_i X = 0$ for $i > n$. Then this holds for all $\PP_k X$, $\PP^k X$, and $\PP^k _k X$.
\end{thm}
\begin{proof}
It suffices to prove this for $\PP^k X$. If $n=-1$ or $0$ then $X \leq n$ and this is trivial, so suppose $n > 0$. We proceed by induction on $|G|$. Examining the connectivities of the slice cells, we have $X \leq n|G|$, so we may assume that $k < n|G|$. Consider the fiber sequence
\begin{align*}
	F \to \PP^k X \to Post^n (\PP^k X).
\end{align*}
Now, by the inductive hypothesis, the map on the right is an isomorphism on all proper subgroups, so $F$ is concentrated over $G$. Furthermore, $F=Post_{n+1} (\PP^k X)$ is $n$-connected, so by Corollary~\ref{cor:GoverGslicefilt} we have
\begin{align*}
	F \geq (n+1)|G| > k.
\end{align*}

\noindent Since $F > k$ and $\PP^k X \leq k$, the map on the left in the above fiber sequence is zero. However, it is also an isomorphism on homotopy groups in dimensions that are $> n$.
\end{proof}

\noindent We now have the following collection of results:
\begin{itemize}
\item A spectrum is zero if and only if its (R)SSS is zero.
\item A spectrum is $n$-connected ($n \leq 0$) if and only if its RSSS is.
\item A spectrum is $n$-coconnected ($n \geq 0$) if and only if its RSSS is.
\item A spectrum restricts to zero in a subgroup of $G$ if and only if its (R)SSS does.
\end{itemize}

\noindent Thus, for example, one will not be forced to compute any of the negative columns in the SS past the first page if they are going to converge to zero anyway. We sum this up by saying that "\emph{the (regular) slice spectral sequence is very efficient}."

\vspace{1 cm}

\indent \emph{Acknowledgement:} I would like to thank Mike Hill for his useful remarks on the present work.


\vspace{1 cm}

\begin{thebibliography}{12}

\bibitem{KR} Daniel Dugger, An Atiyah-Hirzebruch spectral sequence for KR-theory, K-Theory 35 (2005), no. 3-4, 213–256 (2006). MR MR2240234 (2007g:19004)

\bibitem{DF} E. Dror Farjoun. Cellular spaces, null spaces and homotopy localization, volume 1622 of Lecture Notes in Mathematics. Springer-Verlag, New York, 1996.

\bibitem{Primer} Michael A. Hill. The Equivariant Slice Filtration: a Primer. arXiv:1107.3582v1 (2011).

\bibitem{HHR} Michael A. Hill, Michael J. Hopkins, and Douglas C. Ravenel. On the non-existence of elements of Kervaire invariant one. arXiv:0908.3724v2 (2009).

\bibitem{LMS} L. G. Lewis, J. P. May, and M. Steinberger. Equivariant Stable Homotopy Theory, volume 1213 of Lecture Notes in Mathematics. Springer-Verlag, New York, 1986.

\bibitem{Voev} Vladimir Voevodsky. Open problems in the motivic stable homotopy theory. I, Motives, polylogarithms and Hodge theory, Part I (Irvine, CA, 1998), Int. Press Lect. Ser., vol. 3, Int. Press, Somerville, MA, 2002, pp. 3–34. MR MR1977582 (2005e:14030)

\end{thebibliography}
\end{document}